# A PROOF OF "GOLDBACH'S CONJECTURE"

## By Roger Ellman ©

GOLDBACH'S CONJECTURE states:

Every even number greater than two can be expressed as the sum of two primes.

### STEP 1 - General

All of the prime numbers other than *2* are odd. The sum of any two of those odd prime numbers is always an even number. Therefore, it only remains to show that the combinations* of all prime numbers other than *2*, taken two at a time, summed in pairs, yields all of the even numbers greater than *4*. That, along with that the even number *4* is the sum of the pair of prime numbers *[2 + 2]* will complete the proof.

*[Here the combinations in pairs may include the same number twice.]

It can be observed that Goldbach's Conjecture has been verified up to $10^8$ by numerical calculations.[1]

### STEP 2 - How combinations of primes summed in pairs might yield all even numbers.

The odd prime numbers comprise a string of odd numbers each greater than the prior by two except that there are various gaps [intervals of one or more non-primes] in the sequence. For example, all of the odd numbers *< 200* are as follows with the non-primes in bold face [the number *1*, which is divisible only by *1* and itself, is nevertheless defined as non-prime].

[1]
| | | | | | | | | | |
|---|---|---|---|---|---|---|---|---|---|
| **1** | 3 | 5 | 7 | **9** | 11 | 13 | **15** | 17 | 19 | **21** |
| 23 | **25** | **27** | 29 | 31 | **33** | **35** | 37 | **39** | 41 | 43 |
| **45** | 47 | **49** | **51** | 53 | **55** | **57** | 59 | 61 | **63** | **65** |
| 67 | **69** | 71 | 73 | **75** | **77** | 79 | **81** | 83 | **85** | **87** |
| 89 | **91** | **93** | **95** | 97 | **99** | 101 | 103 | **105** | 107 | 109 |
| **111** | 113 | **115** | **117** | **119** | **121** | **123** | **125** | 127 | **129** | 131 |
| **133** | **135** | 137 | 139 | **141** | **143** | **145** | **147** | 149 | 151 | **153** |
| **155** | 157 | **159** | **161** | 163 | **165** | 167 | **169** | **171** | 173 | **175** |
| **177** | 179 | 181 | **183** | **185** | **187** | **189** | 191 | 193 | **195** | 197 |
| 199 | ... | | | | | | | | | |

Designating the set of all prime numbers to be $\{P_i, i = 1, 2, ... \infty\}$ = *3, 5, 7, not 9, 11, ...*, the first sub-set of the set of all combinations* of $P_i$ taken and summed in pairs is the sub-set $\{3 + P_i\}$. That sub-set produces some, but not all, of the even numbers generated as the sum of two primes as follows. Bold face indicates gaps in the even numbers sequence, even numbers that are the sum of two odd numbers one or both of which is non-prime and therefore do not satisfy the conjecture.

[2]   Sub-Set $\{3 + P_i\}$ =

| | | | | | | | | | |
|---|---|---|---|---|---|---|---|---|---|
| 4 | 6 | 8 | 10 | **12** | 14 | 16 | **18** | 20 | 22 | **24** |
| 26 | **28** | 30 | 32 | 34 | **36** | 38 | 40 | **42** | 44 | 46 |
| **48** | 50 | **52** | 54 | 56 | **58** | 60 | 62 | 64 | **66** | 68 |
| 70 | **72** | 74 | 76 | **78** | 80 | 82 | **84** | 86 | 88 | 90 |
| 92 | **94** | **96** | **98** | 100 | **102** | 104 | 106 | **108** | 110 | 112 |
| **114** | 116 | **118** | 120 | 122 | **124** | 126 | **128** | 130 | **132** | 134 |
| **136** | **138** | 140 | 142 | **144** | **146** | 148 | 150 | 152 | 154 | **156** |
| **158** | 160 | **162** | **164** | 166 | **168** | 170 | **172** | **174** | 176 | **178** |
| 180 | 182 | 184 | **186** | **188** | 190 | **192** | 194 | 196 | **198** | 200 |
| ... | | | | | | | | | | |

The gaps in the sequence of even numbers generated by the sub-set $\{3 + P_i\}$ are due to the gaps in the sequence of primes $\{P_i\}$ per equation 1, above. The next sub-set, $\{5 + P_i\}$, fills in some of those gaps while leaving corresponding other ones, again in bold face.

[3]   Sub-Set $\{5 + P_i\}$ =

```
  6    8   10   12   14   16   18   20   22   24   26
 28   30   32   34   36   38   40   42   44   46   48
 50   52   54   56   58   60   62   64   66   68   70
 72   74   76   78   80   82   84   86   88   90   92
 94   96   98  100  102  104  106  108  110  112  114
116  118  120  122  124  126  128  130  132  134  136
138  140  142  144  146  148  150  152  154  156  158
160  162  164  166  168  170  172  174  176  178  180
182  184  186  188  190  192  194  196  198  200  202
204  ...
```

Applying the two sub-sets together, however, the number of gaps in the sequence of all even numbers that are generated as the sum of two primes is reduced. Sub-set $\{5 + P_i\}$ generates some even numbers that are the sum of two primes that $\{3 + P_i\}$ does not. The combined effect is as below, for which if a number is not bold in one or both of equation 2 and 3 then it is not bold below and is not a gap.

[3a]  The even numbers [≤ 200] as generated by Sub-Sets $\{3 + P_i\}$
      and $\{5 + P_i\}$ combined =

```
  4    6    8   10   12   14   16   18   20   22   24
 26   28   30   32   34   36   38   40   42   44   46
 48   50   52   54   56   58   60   62   64   66   68
 70   72   74   76   78   80   82   84   86   88   90
 92   94   96   98  100  102  104  106  108  110  112
114  116  118  120  122  124  126  128  130  132  134
136  138  140  142  144  146  148  150  152  154  156
158  160  162  164  166  168  170  172  174  176  178
180  182  184  186  188  190  192  194  196  198  200
```

The sub-set $\{5 + P_i\}$ supplies a missing even number wherever it encounters the beginning of a gap in the sequence of even numbers that were generated by the $\{3 + P_i\}$ sub-set. That happens because each number in $\{5 + P_i\}$ is 2 more [one odd number higher] than the corresponding number in $\{3 + P_i\}$. The effect is that any number in equation 2 that is at the beginning of a gap [bold face with non-bold face to its left] moves one position to the left in equation 3, moves to a position not in a gap [non-bold face]. The only exception is the initial gap, **4**, which has already been addressed.

Because of this, were the sequence 3, 5, ... of the subsets $\{3 + P_i\}$, $\{5 + P_i\}$, ... to continue without any breaks, for example were it to proceed $\{7 + P_i\}$, $\{\underline{9} + P_i\}$ $\{11 + P_i\}$, ... then all of those sub-sets collectively would eventually fill in all of the gaps in the original sequence generated by $\{3 + P_i\}$ and generate all of the even numbers as sums of two primes. That is, that would happen provided that each of the gaps is such that there are more odd numbers, *(not 3, which generates the original gaps), 5, 7, $\underline{9}$, 11, ...*, preceding the gap than there are even numbers in the gap.

However, the sequence of subsets has breaks such as the $\{\mathbf{9} + P_i\}$ sub-set, which generates no even numbers as the sum of two primes because one of the two numbers summed is always the non-prime **9**. There now remain two issues: what of the invalid sub-sets, ones that produce breaks in the sub-set sequence such as $\{\mathbf{9} + P_i\}$, and are there always sufficient prime [not merely odd] numbers preceding each gap ?

STEP 3 - The invalid sub-sets such as sub-set $\{9 + P_i\}$.

The equations below present the even numbers generated by each of sub-sets $\{7 + P_i\}$ through $\{17 + P_i\}$ for which the index, $I$, where $I = 7$, **9**, 11, 13, 17, ... , is prime and for each presents the cumulative effect of all of the sub-sets through that point in eliminating gaps.

[4]   Sub-Set $\{7 + P_i\}$ =

```
     8    10    12    14    16    18    20    22    24    26    28
    30    32    34    36    38    40    42    44    46    48    50
    52    54    56    58    60    62    64    66    68    70    72
    74    76    78    80    82    84    86    88    90    92    94
    96    98   100   102   104   106   108   110   112   114   116
   118   120   122   124   126   128   130   132   134   136   138
   140   142   144   146   148   150   152   154   156   158   160
   162   164   166   168   170   172   174   176   178   180   182
   184   186   188   190   192   194   196   198   200   ...
```

[4a]  The even numbers [≤ 200] as generated by Sub-Sets $\{3 + P_i\}$, $\{5 + P_i\}$ and $\{7 + P_i\}$ combined =

```
     4     6     8    10    12    14    16    18    20    22    24
    26    28    30    32    34    36    38    40    42    44    46
    48    50    52    54    56    58    60    62    64    66    68
    70    72    74    76    78    80    82    84    86    88    90
    92    94    96    98   100   102   104   106   108   110   112
   114   116   118   120   122   124   126   128   130   132   134
   136   138   140   142   144   146   148   150   152   154   156
   158   160   162   164   166   168   170   172   174   176   178
   180   182   184   186   188   190   192   194   196   198   200
```

[5]   Sub-Set $\{11 + P_i\}$ =

```
    12    14    16    18    20    22    24    26    28    30    32
    34    36    38    40    42    44    46    48    50    52    54
    56    58    60    62    64    66    68    70    72    74    76
    78    80    82    84    86    88    90    92    94    96    98
   100   102   104   106   108   110   112   114   116   118   120
   122   124   126   128   130   132   134   136   138   140   142
   144   146   148   150   152   154   156   158   160   162   164
   166   168   170   172   174   176   178   180   182   184   186
   188   190   192   194   196   198   200   ...
```

[5a]  The even numbers [≤ 200] as generated by Sub-Sets $\{3 + P_i\}$, $\{5 + P_i\}$, $\{7 + P_i\}$ and $\{11 + P_i\}$ combined =

```
     4     6     8    10    12    14    16    18    20    22    24
    26    28    30    32    34    36    38    40    42    44    46
    48    50    52    54    56    58    60    62    64    66    68
    70    72    74    76    78    80    82    84    86    88    90
    92    94    96    98   100   102   104   106   108   110   112
   114   116   118   120   122   124   126   128   130   132   134
   136   138   140   142   144   146   148   150   152   154   156
   158   160   162   164   166   168   170   172   174   176   178
   180   182   184   186   188   190   192   194   196   198   200
```

*[6]*   Sub-Set $\{13 + P_i\}$ =

|  |  |  |  |  |  |  |  |  |  |
|---|---|---|---|---|---|---|---|---|---|
| **14** | 16 | 18 | 20 | **22** | 24 | 26 | **28** | 30 | 32 | 
| | | | | | | | | | **34** |

     **14**   16   18   20   **22**   24   26   **28**   30   32   **34**
     36   **38**   **40**   42   44   **46**   **48**   50   **52**   54   56
     **58**   60   **62**   **64**   66   **68**   70   72   74   **76**   **78**
     80   **82**   84   86   **88**   90   92   **94**   96   **98**   100
     102   **104**   **106**   **108**   110   **112**   114   116   **118**   120   122
     **124**   126   **128**   130   **132**   **134**   **136**   **138**   140   **142**   144
     **146**   **148**   150   152   **154**   **156**   **158**   **160**   162   164   **166**
     **168**   170   **172**   **174**   176   **178**   180   **182**   **184**   186   **188**
     **190**   192   194   **196**   **198**   **200**   ...

*[6a]*  The even numbers [≤ 200] as generated by Sub-Sets $\{3+P_i\}$,
        $\{5+P_i\}$, $\{7+P_i\}$, $\{11+P_i\}$ and $\{13+P_i\}$ combined =

     **4**   6   8   10   12   14   16   18   20   22   24
     26   28   30   32   34   36   38   40   42   44   46
     48   50   52   54   56   58   60   62   64   66   68
     70   72   74   76   78   80   82   84   86   88   90
     92   94   96   **98**   100   102   104   106   108   110   112
     114   116   118   120   122   124   126   **128**   130   132   134
     136   138   140   142   144   146   148   150   152   154   156
     158   160   162   164   166   168   170   172   174   176   178
     180   182   184   186   188   190   192   194   196   198   200

*[7]*   Sub-Set $\{17 + P_i\}$ =

     **18**   20   22   24   **26**   28   30   **32**   34   36   **38**
     40   **42**   **44**   46   48   **50**   **52**   54   **56**   58   60
     **62**   64   **66**   **68**   70   **72**   **74**   76   78   **80**   **82**
     84   **86**   88   90   **92**   **94**   96   **98**   100   **102**   **104**
     106   **108**   **110**   **112**   114   **116**   118   120   **122**   124   126
     **128**   130   **132**   **134**   **136**   **138**   140   **142**   144   **146**   148
     **150**   **152**   154   156   **158**   **160**   **162**   164   166   168   **170**
     **172**   174   **176**   **178**   180   **182**   184   **186**   **188**   190   **192**
     **194**   196   198   **200**   ...

*[7a]*  The even numbers [≤ 200] as generated by Sub-Sets $\{3+P_i\}$,
        $\{5+P_i\}$, $\{7+P_i\}$, $\{11+P_i\}$, $\{13+P_i\}$ and $\{17+P_i\}$ combined =

     **4**   6   8   10   12   14   16   18   20   22   24
     26   28   30   32   34   36   38   40   42   44   46
     48   50   52   54   56   58   60   62   64   66   68
     70   72   74   76   78   80   82   84   86   88   90
     92   94   96   **98**   100   102   104   106   108   110   112
     114   116   118   120   122   124   126   **128**   130   132   134
     136   138   140   142   144   146   148   150   152   154   156
     158   160   162   164   166   168   170   172   174   176   178
     180   182   184   186   188   190   192   194   196   198   200

   Each of the individual single sub-set tables is identical to its predecessor except that in each successive table the number in each position is increased by the amount that the index has increased. The effect is that the numbers in the sequence of tables move continuously to the left and upward in the tables while the table positions of unsatisfactory even numbers [ones in bold face because they are in the positions where there is no number in $P_i$, where the corresponding number in the sequence of all odd numbers is a non-prime] remain unmoved.

The original sequence, that of sub-set *{3 + P<sub>i</sub>}*, sets the problem. It exhibits many even numbers as a sum of two primes, that is as satisfactory numbers, and it also exhibits many gaps, uninterrupted sequences of unsatisfactory even numbers. The length of a gap, *G*, is the number of unsatisfactory numbers in uninterrupted sequence. In general, then, what is required to insure that every even number is generated as the sum of two primes ? What is required to clear all of the gaps ?

As the numbers in the individual single sub-set tables move to the left and upward with the increases in the sub-set index, *I*, whenever a number from a gap moves onto the position of a prime [a non-bold face position] the number is expressed as a sum of two primes, a satisfactory even number, and the gap is reduced to that extent. Not more than *G* such events would be needed to clear the gap. Each valid sub-set that acts after the original *{3 + P<sub>i</sub>}* is a candidate to provide such an event for each gap and whether it succeeds or not it moves the gap nearer to the beginning of the table.

From this point of view, clearing a gap of length *G* requires that the original odd number sequence, equation *1*, exhibit at least *G* primes ahead of the gap. This requirement is conservative because sometimes the same position in a table clears more than one element of a gap, that is one or more of the *G* primes ahead of the gap in the table may sometimes act more than once on the same gap. [For example, for the *G = 3* gap **94 96 98** per sub-set *{3 + P<sub>i</sub>}*, equation *2*, it turns out that the prime *92* actually clears two elements in the gap, the **94** and **96** (see equations *3*, *3a*, *4*, and *4a*). And, for another example, equation *5a* versus equation *6a* where **122** and **126** are cleared simultaneously.] That type of help in clearing gaps is ignored in the present analysis so as to arrive at a worst case.

Another effect may further increase the above requirement. It can happen that a number in a gap could move onto the position of a prime in a sub-set having a non-prime index; or rather, the number in the gap could move more than one number to the left at one time because of the index skipping an intervening invalid sub-set [see equations *4* and *5* between which an invalid sub-set is skipped and likewise equations *6* and *7*]. However viewed, such an event could "waste" a prime, result in failure to benefit from it by the conversion of a number from unsatisfactory to satisfactory, meaning that there may need to be more than *G* available primes preceding the gap. [For example, the **98** in sub-set *{7 + P<sub>i</sub>}*, equation *4*, is in position to move to the left onto a position [non-bold face] that would clear it; however, the next sub-set is *{9 + P<sub>i</sub>}*, an invalid sub-set that is skipped over because its non-prime index prevents clearing any elements. The next position after that is in bold face again and unable to clear an element. The net effect is that the potentially clearing position is missed. "wasted".]

How many more available primes might be needed to account for and offset this effect ? A conservative estimate of the amount that *G* should be increased would be to multiply *G* by the ratio of the number of non-primes preceding the gap to the number of primes preceding the gap, that is increase *G* to *G·[1 + that ratio]*. That is equivalent to assuming that every invalid sub-set wastes a prime. Sometimes the "wasted prime" position is occupied by a bold face number or sometimes the number in position to move there is already cleared so that, either way nothing is really lost. [For example, the **142** in sub-set *{7 + P<sub>i</sub>}*, equation *4*.] An adjustment to take account of that help in clearing gaps is ignored in the present analysis so as to arrive at a worst case.

STEP *4* - The number of primes required to clear a gap

The revised value of the required number of available primes preceding the gap now becomes as follows.

*[8]*   Where:

   N ≡ the first number of a gap, the gap's beginning.

$\pi(N) \equiv$ the number of primes that are $\leq N$.

$R(N) \equiv$ the required number of primes ahead of the gap that begins at number "N" needed to clear that gap.

Then:

$$R(N) = G \cdot \left[1 + [N-\pi(N)]/\pi(N)\right]$$

$$= G \cdot \left[N/\pi(N)\right]$$

*STEP 5* - <u>Why there are always sufficient primes preceding each gap.</u>

The subject of the distribution of primes has been studied in depth. Chapter 4, "How are the Prime Numbers Distributed ?" of Reference [1] summarizes the history and results of those studies and presents a number of related proofs. The results fall into two categories for the present purposes: section I of the chapter, which treats the development of the Prime Number Theorem, that is expressions for the number of primes in designated intervals, and section II of the chapter, which treats gaps between primes. The source of the following data is that work. [1]

The Prime Number Theorem, the most fundamental theorem of prime numbers is as follows.

[9]  $\pi(N) \equiv$ the prime counting function

$\equiv$ the number of primes $\leq N$, an integer

$\approx N/Ln(N)$ the approximation improving as $N$ increases

That approximation is low by *5.78%* for $N = 10^8$ improving to low by *2.79%* for $N = 10^{16}$. A better approximation is given by a function called the logarithmic interval as follows.

[10]  $\pi(N) \approx Li(N) = \int_2^N dx/Log(x)$

The logarithmic interval approximation to $\pi(N)$ is high by only *0.013%* for $N = 10^8$ and by only *0.000,000,5%* for $N = 10^{16}$.

Even more accurate is the Riemann function, too involved to be worth specifying here, which for $N = 10^8$ differs from the correct value by only *0.0017%* and for $N = 10^{16}$ by *0.000,000,1%*.

Substituting equation *9* into equation *8* the following is obtained.

[11]  $R(N) = G \cdot [N/\pi(N)]$

$= G \cdot N \cdot [Ln(N)/N]$

$= G \cdot Ln(N)$

The issue is, of course, how does *R(N)*, the number of preceding primes required, compare with $\pi(N)$, the number of preceding primes actually available ? That is as follows.

[12]

$$\pi(N)/R(N) = \frac{N/Ln(N)}{G \cdot Ln(N)}$$

$$= N/G \cdot [Ln(N)]^2$$

Figure 1, below lists some values for that function using values for *G* extrapolated from the percent deviation of the logarithmic interval from the exact count of *π(N)*

presented above. From hundreds to thousands to even far greater multiples of the required number of primes preceding the gaps are actually available for clearing them.

| N | G | $\pi(N)/R(N)$ | $G/N$ % |
|---|---|---|---|
| $10^8$ | $10^{3.0}$ | 295 | 0.001 |
| $10^9$ | $10^{3.6}$ | 585 | 0.0004 |
| $10^{10}$ | $10^{4.3}$ | 945 | 0.0002 |
| $10^{11}$ | $10^{4.9}$ | 1962 | 0.00008 |
| $10^{12}$ | $10^{5.5}$ | 4142 | 0.00003 |
| $10^{13}$ | $10^{6.1}$ | 8865 | 0.00001 |
| $10^{14}$ | $10^{6.8}$ | 15251 | 0.000006 |
| $10^{15}$ | $10^{7.4}$ | 33372 | 0.000003 |
| $10^{16}$ | $10^{8.0}$ | 73676 | 0.000001 |

*Figure 1*

The formulations of various accuracies for evaluating $\pi(N)$ cited on the previous page: $N/Ln(N)$, the logarithmic interval and the Riemann function, cannot be used to find individual prime numbers. Precise, not limited accuracy, is required for that. However, the formulations do indicate that the prime numbers are somewhat well distributed over the range of all numbers, that dense concentrations and large gaps cannot occur. The function of equation *9* is smooth as shown in Figure 2, below, and while the function only approximates the actual values of $\pi(N)$, that approximation is fairly close over most of the range. The same can be said, and more emphatically, of the more accurate functions cited for which the approximation is better. Gaps that are large relative to the location in the sequence of numbers, that is other than small values of $G/\pi(N)$ simply cannot occur. The precise condition, $\pi(N) \geq R(N) = G \cdot Ln(N)$, is greatly exceeded.

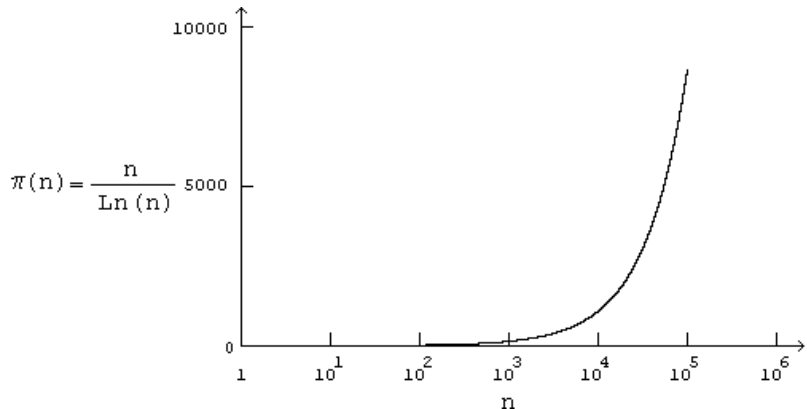

*Figure 2, $\pi(N)$*

For example, at $N = 10^8$, the deviation of $\pi(N)$ from the exact number of primes up to that *N* is *0.0017%*, which is a number range of *1,700* out of *100,000,000*. The largest gap in that neighborhood could not reasonably exceed twice that *3,400*, and certainly not ten times it. The number of primes available up to that point is *5,761,455*. That is far more than needed to eliminate the effect of the gap. [The largest gap that has ever been found is a sequence of *653* non-primes following the prime *11,000,001,446,613,353* at which the number of preceding primes available is more than *279,238,341,033,925* (the value for $10^{16}$).]

*IN SUMMATION:*

1 - All of the prime numbers other than *2* are odd, *2* being the only even prime number.  Further, the even number *4 = 2 + 2*.

2 - The sum of any two of the odd prime numbers is always an even number.

3 - All combinations* of the odd numbers ≥ *3* [whether prime or not] summed in pairs produces all of the even numbers ≥ *6*.

3 - While just the prime odd numbers in sequence is a sequence with gaps as compared to that of all of the odd numbers; nevertheless, all combinations of the odd prime numbers ≥ *3* summed in pairs produces all of the even numbers provided that there are enough primes preceding the gaps.

4 - That requirement is that *π(N) ≥ R(N) = G·Ln(N)* where *N* is the first number in the gap, *π(N)* is the number of primes less than or equal to *N*, *R(N)* is the number of preceding primes needed to assure clearance of the gap, and *G* is the number of sequential non-primes in the gap.  This requirement is comprehensively satisfied by all of the prime numbers and gaps because of the sufficiently smooth nature of *π(N)*.

Which proves the conjecture.

Reference:

[1] Ribenboim, P., *The Book of Prime Number Records, Second Edition*, Springer-Verlag, 1989, Library of Congress catalog # 89-21675.